\documentclass[11pt, a4paper]{article}

\usepackage{amsthm}
\usepackage{amsmath}
\usepackage{amssymb}
\usepackage{amsfonts}
\usepackage{latexsym}
\usepackage{booktabs}
\usepackage{multirow}
\usepackage{graphicx}
\usepackage{natbib}
\usepackage{pdfpages}
\usepackage{algpseudocode}
\usepackage[left = 2.5cm, right = 2.5cm, top = 2.5cm, bottom = 3cm]{geometry}
\usepackage[colorlinks,citecolor=black,urlcolor=black]{hyperref}
\usepackage{xr}
\def\diag{\mathop{\rm diag}}

\newtheoremstyle{example}
{3pt} 
{3pt} 
{} 
{0\parindent} 
{\bf}
{:} 
{.5em} 
{} 
\newtheoremstyle{theorem}
{3pt} 
{3pt} 
{\em} 
{0\parindent} 
{\bf}
{:} 
{.5em} 
{} 
\theoremstyle{example} \newtheorem{example}{Example}
\theoremstyle{theorem} \newtheorem{theorem}{Theorem}
\theoremstyle{lemma} 
\theoremstyle{corollary} \newtheorem{corollary}{Corollary}

\author{Ioannis Kosmidis \\ \texttt{ioannis.kosmidis@warwick.ac.uk} \smallskip
  \\ and \smallskip \\
  David Firth \\ \texttt{d.firth@warwick.ac.uk} \bigskip \\
  Department of Statistics, University of Warwick\\
  Coventry CV4 7AL, UK \smallskip\\
  and \smallskip \\
  The Alan Turing Institute \\
  British Library, London NW1 2DB, UK \\
}

\title{Jeffreys-prior penalty, finiteness and shrinkage in \\ binomial-response
  generalized linear models}

\begin{document}

\maketitle

\begin{abstract}
  Penalization of the likelihood by Jeffreys' invariant prior, or by a
  positive power thereof, is shown to produce finite-valued maximum
  penalized likelihood estimates in a broad class of binomial
  generalized linear models.  The class of models includes logistic
  regression, where the Jeffreys-prior penalty is known additionally
  to reduce the asymptotic bias of the maximum likelihood estimator;
  and also models with other commonly used link functions such as
  probit and log-log.  Shrinkage towards equiprobability across
  observations, relative to the maximum likelihood estimator, is
  established theoretically and is studied through illustrative
  examples.  Some implications of finiteness and shrinkage for
  inference are discussed, particularly when inference is based on
  Wald-type procedures.  A widely applicable procedure is developed
  for computation of maximum penalized likelihood estimates, by using
  repeated maximum likelihood fits with iteratively adjusted binomial
  responses and totals.  These theoretical results and methods
  underpin the increasingly widespread use of reduced-bias and
  similarly penalized binomial regression models in many applied
  fields. \\ \noindent {Keywords:
    \textit{logit}; \textit{probit}; \textit{bias reduction};
    \textit{penalized likelihood}; \textit{data separation};
    \textit{infinite estimate}; \textit{working weight};
    \textit{Bradley-Terry model}}
\end{abstract}

\section{Introduction}

Logistic regression is one of the most frequently applied generalized
linear models in statistical practice, both for inference about
covariate effects on binomial probabilities, and for prediction.
Consider realizations $y_1, \ldots, y_n$ of independent binomial
random variables $Y_1, \ldots, Y_n$ with success probabilities
$\pi_1, \ldots, \pi_n$ and totals $m_1, \ldots, m_n$, respectively.
Suppose that each $y_i$ is accompanied by a $p$-dimensional covariate
vector $x_i$ and that the model matrix $X$ with rows
$x_1, \ldots, x_n$ has full rank. A logistic regression model has
\begin{equation}
  \label{binomialGLM}
  \pi_i = (G \circ \eta_i)(\beta) \quad \text{with} \quad G(\eta) = \frac{e^\eta}{1 + e^\eta} \quad \text{and} \quad \eta_i(\beta) = \sum_{t = 1}^p \beta_tx_{it} \quad (i = 1, \ldots, n)\, ,
\end{equation}
where $\beta = (\beta_1, \ldots, \beta_p)^\top$ is the $p$-dimensional
parameter vector, and $x_{it}$ is the $t$th element of $x_i$
$(i = 1, \ldots, n)$; if an intercept parameter is present in the
model then the first column of $X$ is a vector of ones.
The maximum likelihood estimator $\hat\beta$ of $\beta$ in
(\ref{binomialGLM}) maximizes the
log-likelihood
\begin{equation}
\label{loglik}
l(\beta) = \sum_{i = 1}^n y_i \eta_i(\beta) - \sum_{i = 1}^n m_i \log\left\{1 + e^{\eta_i(\beta)}\right\}\, .
\end{equation}

\citet{albert+anderson:1984} categorized the possible settings for the
sample points $(y_1, x_1^\top)^\top$, $\ldots$, $(y_n, x_n^\top)^\top$
into complete separation, quasi-complete separation and
overlap. Specifically, if there exists a vector $\gamma \in \Re^p$
such that $\gamma^\top x_i > 0$ for all $i$ with $y_i > 0$ and
$\gamma^\top x_i < 0$ for all $i$ with $y_i = 0$, then there is
complete separation in the sample points; if there exists a vector
$\gamma \in \Re^p$ such that $\gamma^\top x_i \ge 0$ for all $i$ with
$y_i > 0$ and $\gamma^\top x_i \le 0$ for all $i$ with $y_i = 0$, then
there is quasi-complete separation in the sample points; and if
neither complete nor quasi-complete separation is present, then the
sample points overlap.  \citet{albert+anderson:1984} showed that
separation is necessary and sufficient for the maximum likelihood
estimate to have at least one infinite-valued component. A parallel
result appears in \citet{silvapulle:1981}, where it is shown that if
$G(\eta)$ in~(\ref{binomialGLM}) is any strictly increasing
distribution function such that $-\log G(\eta)$ and $\log\{1- G(\eta)\}$
are convex, and $x_{i1} = 1$ $(i = 1, \ldots, n)$, then the maximum
likelihood estimate has all components finite if and only if there is
overlap.

When data separation occurs, standard maximum-likelihood estimation
procedures, such as iteratively reweighted least squares
\citep{green:1984}, can be numerically unstable
due to the occurrence of large parameter values as
the procedures attempt to maximize~(\ref{loglik}). In addition,
inferential procedures that directly depend on the estimates and the
estimated standard errors, such as Wald tests, can give misleading
results. For a recent review of such problems and some solutions, see
\citet{mansournia+geroldinger+greenland+heinze:2018}.

\citet{firth:1993} showed that if the logistic
regression likelihood is penalized by Jeffreys' invariant prior, then the resulting
maximum penalized likelihood estimator has bias of smaller
asymptotic order than that of the maximum likelihood estimator in
general. Specifically, for logistic regressions the reduced-bias
estimator $\tilde\beta$ results from maximization of
\begin{equation}
  \label{penloglik}
  \tilde{l}(\beta) = l(\beta) + \frac{1}{2}\log\left|X^\top W(\beta)X\right| \, ,
\end{equation}
with $W(\beta) = \diag\{w_1(\beta), \ldots, w_n(\beta)\}$, and where
$w_i(\beta) = m_i(\omega \circ \eta_i)(\beta)$ is the working weight for
the $i$th observation with
$\omega(\eta) = e^\eta/(1 + e^\eta)^2$ $(i = 1, \ldots, n)$.
\citet{heinze+schemper:2002}, in extensive numerical studies, observed
that the reduced-bias estimates have finite values even when data
separation occurs. Based on an argument about parameter-dependent
adjustments to $y_1, \ldots, y_n$ and $m_1, \ldots, m_n$ stemming from
the form of the gradient of~(\ref{penloglik}),
\citet{heinze+schemper:2002} conjectured that finiteness of the
reduced-bias estimates holds for every combination of data and
logistic regression model. \citet{heinze+schemper:2002} also observed
that the reduced-bias estimates are typically smaller in absolute
value than the corresponding maximum likelihood estimates, when the
latter are finite. These observations are in agreement with the
asymptotic bias of the maximum likelihood estimator in logistic
regressions being approximately collinear with the parameter vector
\citep[see, for example,][]{cordeiro+mccullagh:1991}.

Example~\ref{bt} illustrates the finiteness and shrinkage properties
of the maximum penalized likelihood estimator in the context of
estimating the strength of NBA basketball teams using a Bradley-Terry
model \citep{bradley+terry:1952}.

\begin{example}
  \label{bt}

  Suppose that $y_{ij} = 1$ when team $i$ beats team $j$, and
  $y_{ij} = 0$, otherwise. The Bradley-Terry model assumes that the
  contest outcome $y_{ij}$ is the realization of a Bernoulli random
  variable with probability
  $\pi_{ij} = \exp(\beta_i - \beta_j)/\{1 + \exp(\beta_i -
  \beta_j)\}$, and that the outcomes for the available contests are
  independent. The Bradley-Terry model is a logistic regression with
  probabilities as in~(\ref{binomialGLM}), for the particular $X$
  matrix whose rows are indexed by contest identifiers $(i, j)$ and
  whose general element is
  $x_{ij,t} = \delta_{it} - \delta_{jt}\quad (t = 1, \ldots, p)$.
  Here, $\delta_{it}$ is the Kronecker delta, with value one when
  $t = i$ and zero otherwise. The parameter $\beta_t$ can be thought
  as measuring the ability or strength of team $t$
  $(t = 1, \ldots, p)$. Only contrasts are estimable, and an
  identifiable parameterization can be achieved by setting one of the
  abilities to zero.  See, for example, \citet[\S~11.6]{agresti:2013}
  for a general discussion of the model.

  We use the Bradley-Terry model to estimate the ability of basketball
  teams from game outcomes in the regular season of the 2014--2015 NBA
  conference. For illustrative purposes, we use only the $262$ games
  that took place before 3 December 2014, up to which date the
  Philadelphia 76ers had recorded 17 straight losses and no win.  The
  dataset was obtained from \texttt{www.basketball-reference.com} and is
  also available as part of the Supplementary Material. The ability of
  the San Antonio Spurs, the champion team of the 2013--2014
  conference, is set to zero, so that each $\beta_i$ is the contrast
  of the ability of team $i$ with San Antonio Spurs. The model is
  estimated via iteratively reweighted least squares, as implemented
  in the \texttt{glm} function of R \citep{rproject} with default
  settings for the optimization. No warning or error was returned
  during the fitting process.

\begin{figure}[t!]
  \begin{center}
    \includegraphics[width = \textwidth]{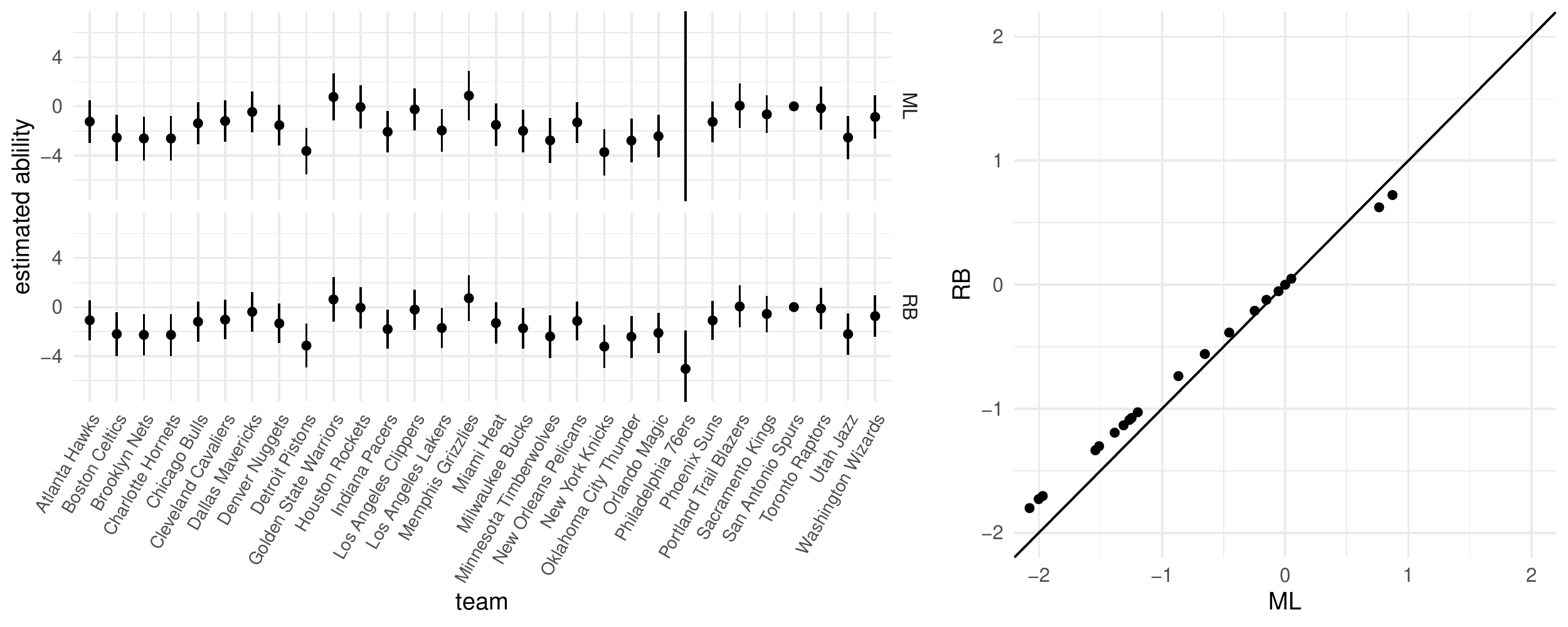}
  \end{center}
  \caption{Left: Estimated contrasts with San Antonio Spurs, in
    ability of NBA teams. The abilities are estimated using a
    Bradley-Terry model on the outcomes of the 262 games before 03
    December 2014 in the regular season of the 2014--2015 NBA
    conference, using the maximum likelihood (ML, top) and the
    reduced-bias (RB, bottom) estimator. The vertical segments are
    nominally 95\% Wald-type confidence intervals. Right: reduced-bias
    estimates of ability contrasts versus maximum likelihood estimates
    of ability contrasts. The maximum likelihood estimate for
    Philadelphia 76ers is not plotted. The solid $45^o$ line is
    displayed for reference.  }
  \label{abilities}
\end{figure}
  
  The top panel in the left plot of Figure~\ref{abilities} shows
  the reported maximum likelihood estimates of the contrasts, along with their
  corresponding nominally $95\%$ individual Wald-type confidence
  intervals. The contrast for Philadelphia 76ers stands out
  in the output from \texttt{glm} with a value of $-19.24$ and a
  corresponding estimated standard error of $844.97$. These values are
  in fact representations of $-\infty$ and $\infty$, respectively, as
  confirmed by the \texttt{detect\_separation} method of the
  \texttt{brglm2} R package \citep{brglm2}, which implements
  separation-detection algorithms from a 2007 University of Oxford
  Department of Statistics PhD thesis by K.~Konis.
  The data are separated, with the maximum
  likelihood estimates for all teams being finite except that for
  Philadelphia 76ers, which is minus infinity. A particularly worrying
  side-effect of data separation here is that if the computer output
  is used naively, a Wald test for difference in ability between
  Philadelphia 76ers and San Antonio Spurs results in no apparent
  evidence of a difference, which is counter intuitive given that the
  former had no wins in 17 games and the latter had 13 wins in 17
  games.  In contrast, the reduced-bias estimates in the bottom panel
  of the left of Figure~\ref{abilities} all have finite values and
  finite standard errors. The right plot in~Figure~\ref{abilities}
  illustrates the shrinkage of the reduced-bias estimates towards zero
  that has also been discussed in a range of different settings, for
  example in \citet{heinze+schemper:2002} and \citet{zorn:2005}.

\end{example}

The apparent finiteness and shrinkage properties of the reduced-bias
estimator, coupled with the fact that the estimator has the same
first-order asymptotic distribution as the maximum likelihood
estimator, are key reasons for the increasingly widespread use of
Jeffreys-prior penalized logistic regression in applied work.  At the
time of writing, \emph{Google Scholar} records approximately 2700
citations of \citet{firth:1993}, more than half of which are from 2015
or later.  The list of application areas is diverse, including for
example agriculture and fisheries research, animal and plant ecology,
criminology, commerce, economics, psychology, health and medical
sciences, politics and many more.  The particularly strong uptake of
the method in health and medical sciences, and politics stems largely
from the works of \citet{heinze+schemper:2002} and \citet{zorn:2005},
respectively. The reduced-bias estimator is also implemented in
dedicated open-source software, such as the \texttt{brglm2}
\citep{brglm2} and \texttt{logistf} \citep{logistf} R packages, and it
has now become part of textbook treatments of logistic regression;
see, for example, \citet[\S~7.4]{agresti:2013}, or
\citet[\S~10.3]{hosmer+lemeshow:2013}.

However, a definitive theoretical account of the empirically evident
finiteness and shrinkage properties has yet to appear in the
literature. Such a formal account is much needed, particularly in
light of recent advances that demonstrate benefits of the reduced-bias
estimator in wider contexts than the ones for which it was originally
developed. An example of such an advance is \citet{lunardon:2018},
which explores the performance of bias reduction in stratified
settings and shows that bias reduction is particularly effective for
inference about a low-dimensional parameter of interest in the
presence of high-dimensional nuisance parameters. For the estimation
of high-dimensional logistic regression models with $p/n \to \kappa$,
$\kappa \in (0, 1)$, experiments reported in the supplementary
information of \citet{sur+candes:2019} (see, also,
Section~S3.3 in the Supplementary Material) show
that bias reduction performs similarly to their newly proposed method,
and markedly better than maximum likelihood. These new theoretical and
empirical results justify and motivate use of the reduced-bias
estimator in even more complex applied settings than the one covered
by the framework of \citet{firth:1993}; in such settings, more
involved methods such as modified profile likelihoods \citep[see, for
example][]{sartori:2003} and approximate message-passing algorithms
\citep[see, for example][]{sur+candes:2019} have also been proposed
for recovering inferential accuracy.

This paper formally derives the finiteness and shrinkage properties of
reduced-bias estimators for logistic regressions under only the
condition that model matrix $X$ has full rank. We also provide
geometric insights on how penalized likelihood estimators shrink
towards zero, and discuss the implications of finiteness and shrinkage
in inference, especially in hypothesis tests and confidence regions
using Wald-type procedures.

It is shown how the results extend in a direct way to other
commonly-used link functions, such as the probit, log-log,
complementary log-log and cauchit, whenever the Jeffreys prior is used
as a likelihood penalty. The work presented here thus complements
earlier work of \citet{ibrahim+laud:1991} and especially
\citet{chen+ibrahim+kim:2008}, which studies the same models from a
Bayesian perspective. Here we study the behaviour of the posterior
mode and thereby derive results that add to those earlier findings,
whose focus was instead on important Bayesian aspects such as
propriety and moments of the posterior distribution.

The results in this paper also readily extend to situations where
penalized log-likelihoods of the form
\begin{equation}
  \label{penloglika}
  l^\dagger(\beta; a) = l(\beta) + a \log\left| X^\top W(\beta)X\right| \quad (a > 0)\, ,
\end{equation}
are used, with $a$ allowed to take values other than $1/2$. Such
penalized log-likelihoods have proved useful in prediction contexts,
where the value of $a$ can be tuned to deliver better estimates of the
binomial probabilities; and they are the subject of ongoing research
\citep[see, for example,][]{elgmati+fiaccone+henderson+matthews:2015,
  puhr+heinze+nold+lusa+geroldinger:2017}.  The repeated maximum
likelihood fits procedure with iteratively adjusted binomial responses
and totals, derived in Section~\ref{algorithm}, maximizes
$l^\dagger(\beta; a)$ for general binomial-response generalized linear
models and any $a>0$.

\section{Logistic regression}
\label{logistic}

\subsection{Preamble}

Results on finiteness and shrinkage of the maximum penalized
likelihood estimator are derived first for logistic regression, which
is the leading case in applications and also the case for which
maximum penalized likelihood, with Jeffreys-prior penalty, coincides
with asymptotic bias reduction.  These results provide a platform for
the generalization to link functions other than logit in
Section~\ref{extensions}.

\subsection{Finiteness}
\label{finite}
Let $W^*(r)$ be $W(\beta)$ at $\beta = \beta(r)$, $r\in \Re$, where
$\beta(r)$ is a path in $\Re^p$ such that $\beta(r) \to \beta_0$ as
$r \to \infty$, with $\beta_0$ having at least one infinite
component. Theorem~\ref{vanish} below describes the limiting behaviour
of the determinant of the expected information matrix
$X^\top W^*(r) X$ as $r$ diverges to infinity, only under the
assumption that $X$ is of full rank. An important implication of
Theorem~\ref{vanish} is Corollary~\ref{finiteness} which shows that
the reduced-bias estimators for logistic regressions are always
finite. These new results formalize a sketch argument made in
\citet[\S~3.3]{firth:1993}.

\begin{theorem}
\label{vanish}
Suppose that $X$ has full rank. Then
$\lim_{r \to \infty} |X^\top W^*(r)X| = 0$.
\end{theorem}

\begin{corollary}
  \label{finiteness}
  Suppose that $X$ has full rank. The vector $\tilde\beta$ that
  maximizes $\tilde{l}(\beta)$ has all of its components finite.
\end{corollary}
The proofs of Theorem~\ref{vanish} and Corollary~\ref{finiteness} are
given in the Supplementary Material.

Corollary~\ref{finiteness} also holds for any fixed $a > 0$
in~(\ref{penloglika}). As a result, the maximum penalized likelihood
estimators from the maximization of $l^\dagger(\beta; a)$ in
(\ref{penloglika}) have finite components, for any $a > 0$.

Despite its practical utility, the finiteness of the reduced-bias
estimator results in some notable, and perhaps undesirable,
side-effects on Wald-type inferences based on the reduced-bias
estimator that have been largely overlooked in the literature. The
finiteness of $\tilde\beta$ implies that the estimated standard errors
$s_t(\tilde{\beta})$ $(t = 1, \ldots, p)$, calculated as the square
roots of the diagonal elements of the inverse of
$X^\top W(\tilde\beta)X$, are also always finite. Since
$y_1, \ldots, y_n$ are realizations of binomial random variables,
there is only a finite number of values that the estimator
$\tilde\beta$ can take for any given $x_1, \ldots, x_n$. Hence, there
will always be a parameter vector with large enough components that
the usual Wald-type confidence intervals
$\tilde\beta_t \pm z_{1 - \alpha/2} s_t(\tilde\beta)$, or confidence
regions in general, will fail to cover regardless of the nominal level
$\alpha$ that is used. This has also been observed in the complete
enumerations of \citet{kosmidis:2014} for proportional odds models
which are extensions of logistic regression to ordinal responses; and
it is also true when the penalized likelihood is profiled for the
construction of confidence intervals, as is proposed, for example, in
\citet{heinze+schemper:2002}, and in \citet{bull+lewinger+lee:2007}
for multinomial regression models.

\begin{table}[t]
  \caption{Common link functions and the corresponding forms for
    $G(\eta)$ and $\omega(\eta)$. For all the displayed link
    functions, $\omega(\eta)$ vanishes as $\eta$ diverges.}
\centering
\begin{tabular}{ccc}
  \smallskip
Link function &  $G(\eta)$ & $\omega(\eta)$ \\ \smallskip
logit & $\displaystyle \frac{e^\eta}{1 + e^\eta}$ & $\displaystyle \frac{e^\eta}{(1 + e^\eta)^2}$ \\ \smallskip
probit & $\displaystyle \Phi(\eta)$ & $\displaystyle \frac{\{\phi(\eta)\}^2}{\Phi(\eta)\{1- \Phi(\eta)\}}$ \\ \smallskip
c-log-log & $\displaystyle 1 - e^{-e^\eta}$ & $\displaystyle \frac{e^{2 \eta}}{e^{e^{\eta}}-1}$ \\ \smallskip
log-log & $\displaystyle e^{-e^{- \eta}}$ & $\displaystyle \frac{e^{-2 \eta}}{e^{e^{-\eta}}-1}$ \\ \smallskip
cauchit & $\displaystyle \frac{1}{2} + \frac{\arctan(\eta)}{\pi}$ & $\displaystyle \frac{1}{\left(1 + \eta^2\right)^2  \left[\frac{\pi^2}{4} -  \{\tan^{-1}(\eta)\}^2 \right]}$ \\ 
\end{tabular}
\label{weightsLink}
\end{table}

\subsection{Shrinkage}
\label{shrinkage}
The following theorem is key when exploring the shrinkage properties
of the reduced-bias estimator that have been illustrated in
Example~\ref{bt}.

\begin{theorem}
  \label{JeffreysMax}
  Suppose that $X$ has full rank. Then
  \begin{itemize}
  \item[i)] The function $|X^\top W(\beta)X|$ is globally maximized at $\beta = 0$.
  \item[ii)] If $\bar{W}(\pi) =
    \diag\{m_1\pi_1(1-\pi_1), \ldots, m_n\pi_n(1-\pi_n)\}$, then $|X^\top \bar{W}(\pi)X|$ is log-concave on $\pi$.
  \end{itemize}
\end{theorem}
A complete proof of Theorem~\ref{JeffreysMax} is in the Supplementary
Material. Part i) also follows directly from
\citet[Theorem~1]{chen+ibrahim+kim:2008}.

Consider estimation by maximization of the penalized log-likelihood
$l^\dagger(\beta; a)$ in~(\ref{penloglika}) for $a = a_1$ and
$a = a_2$ with $a_1 > a_2 \ge 0$.  Let $\beta^{(a_1)}$ and
$ \beta^{(a_2)}$ be the maximizers of $l^\dagger(\beta; a_1)$ and
$l^\dagger(\beta; a_2)$, respectively and $ \pi^{(a_1)}$ and
$ \pi^{(a_2)}$ the corresponding estimated $n$-vectors of
probabilities.  Then, by the concavity of
$\log|X^\top \bar{W}( \pi)X|$, the vector $ \pi^{(a_1)}$ is closer to
$(1/2, \ldots, 1/2)^\top$ than is $ \pi^{(a_2)}$, in the sense that
$ \pi^{(a_1)}$ lies within the hull of that convex contour of
$\log |X^\top\bar{W}( \pi)X|$ containing $ \pi^{(a_2)}$. With the
specific values $a_1=1/2$ and $a_2=0$ the last result refers to
maximization of the likelihood penalized by Jeffreys prior
and to maximization of the un-penalized likelihood,
respectively. Hence, use of reduced-bias estimators for logistic
regressions has the effect of shrinking towards the model that implies
equiprobability across observations, relative to maximum
likelihood. Shrinkage here is according to a metric based on the
expected information matrix rather than to Euclidean distance. Hence,
the reduced-bias estimates are only typically, rather than always,
smaller in absolute value than the corresponding maximum likelihood
estimates.

If the determinant of the inverse of the expected information matrix
is considered as a generalized measure of the asymptotic variance,
then the estimated generalized asymptotic variance at the reduced-bias
estimates is always smaller than the corresponding estimated variance
at the maximum likelihood estimates.
Hence approximate confidence ellipsoids, based on asymptotic normality of the reduced-bias estimator, are reduced in volume.

\section{Non-logistic link functions}
\label{extensions}

\subsection{Preamble}

The results here generalize Sections~\ref{finite} and \ref{shrinkage}
beyond the logit link, still for estimators from penalized likelihoods
of form (\ref{penloglika}).
For non-logistic links, such estimators no longer coincide with 
the bias-reduced estimator of \citet{firth:1993}.

\subsection{Finiteness}

The results in Theorem~\ref{vanish} and Corollary~\ref{finiteness}
readily extend to more link functions than the logistic. Specifically,
if $G(\eta) = e^\eta/(1 + e^\eta)$ in model~(\ref{binomialGLM}) is
replaced by an at least twice differentiable and invertible function
$G: \Re \to (0, 1)$, then the expected information matrix has again
the form $X^\top W(\beta)X$ but with working weights
$w_i(\beta) = m_i(\omega \circ \eta_i)(\beta)$ $(i = 1, \ldots, n)$ where
$\omega(\eta) = g(\eta)^2/[G(\eta)\{1 - G(\eta)\}]$ and
$g(\eta) = d G(\eta)/d\eta$. If the link function is such that
$\omega(\eta) \to 0$ as $\eta$ diverges to either $-\infty$ or
$\infty$, then the proofs of Theorem~\ref{vanish} and
Corollary~\ref{finiteness} in the Supplementary Material
apply unaltered to show that
$\lim_{r \to \infty} |X^\top W^*(r)X| = 0$ and, when the penalty
is a positive power of Jeffreys' invariant prior, the maximum
penalized likelihood estimates have finite components. The logit,
probit, complementary log-log, log-log and cauchit links are some
commonly-used link functions for which $\omega(\eta) \to 0$. The
functions $G(\eta)$ and $\omega(\eta)$ for each of the above link
functions are shown in Table~\ref{weightsLink}.

\subsection{Shrinkage}

Let $\bar{\omega}(z) = \{(g \circ G^{-1})(z)\}^2/\{z(1 - z)\}$. If the link
function is such that $\bar{\omega}(z)$ is maximized at some value
$z_0 \in (0, 1)$, then the same arguments as in the proof of result i)
in Theorem~\ref{JeffreysMax} can be used to show that
$|X^\top \bar{W}(\pi)X|$ is globally maximized at
$(z_0, \ldots, z_0)^\top$. The left plot of
Figure~\ref{shrinkage.plot} illustrates that this condition is
satisfied for the logit, probit, log-log, and complementary log-log
link functions. If $x_{i1} = 1$ $(i = 1, \ldots, n)$, then the maximum
of $|X^\top W(\beta)X|$ is achieved at
$\beta = (b_0, 0, \ldots, 0)^\top$, where $b_0 = g^{-1}(z_0)$\null. In addition, directly from
the proof of Theorem~\ref{JeffreysMax}, a sufficient condition for the
log-concavity of $|X^\top \bar{W}(\pi)X|$ for non-logit link functions
is that $\bar{\omega}(z)$ is concave.

\section{Maximum penalized likelihood as repeated maximum likelihood}
\label{algorithm}

The maximum penalized likelihood estimates, for full rank $X$, can be
computed by direct numerical optimization of the penalized
log-likelihood $l^\dagger(\beta; a)$ in~(\ref{penloglika}) or by using
a quasi Newton-Raphson iteration as in
\citet{kosmidis+firth:2010}. Nevertheless, the particular form of the
Jeffreys prior allows the convenient computation of penalized
likelihood estimates by leveraging readily available
maximum-likelihood implementations for binomial-response generalized
linear models.

If $G(\eta) = e^\eta / (1 + e^\eta)$ in model~(\ref{binomialGLM}) is
replaced by any invertible function $G: \Re \rightarrow (0, 1)$ that
is at least twice differentiable, then differentiation of
$l^\dagger(\beta; a)$ with respect to $\beta_t$ $(t = 1, \ldots, q)$
gives that the penalized likelihood estimates are the solutions of
\begin{equation}
  \label{scores_jeffreys} \sum_{i = 1}^n \frac{w_i(\beta)}{d_i(\beta)} \left[ y_i + 2 a h_i(\beta) \left\{ q_i(\beta) - \frac{1}{2} \right\} - m_i \pi_i(\beta) \right] x_{it} = 0 \quad (t = 1, \ldots, p) \, ,
\end{equation}
where $\pi_i(\beta) = (G \circ \eta_i)(\beta)$,
$d_i(\beta) = m_i (g \circ \eta_i)(\beta)$,
$q_i(\beta) = d_i'(\beta)/w_i(\beta) + \pi_i(\beta)$, and
$d_i'(\beta) = m_i (g' \circ \eta_i)(\beta)$ with
$g'(\eta) = d^2 G(\eta) / d\eta^2$. The quantity $h_i(\beta)$
$(i = 1, \ldots, n)$ is the $i$th diagonal element of the `hat' matrix
$H(\beta) = X \{X^\top W(\beta) X\}^{-1} X^\top W(\beta)$.

If we temporarily omit the observation index and suppress the
dependence of the various quantities on $\beta$, the derivatives of
$l^\dagger(\beta; a)$ are the derivatives of the binomial
log-likelihood $l(\beta)$ with link function $G(\eta)$, after
adjusting the binomial response $y$ to $y + 2 a h (q - 1/2)$. Hence,
the penalized likelihood estimates can be conveniently computed
through repeated maximum-likelihood fits, where each repetition
consists of two steps: P1) the adjusted responses are computed at the
current parameter values; and P2) the maximum likelihood estimates of
$\beta$ are computed at the current value of the adjusted responses.

However, depending on the sign and magnitude of $2 a h (q - 1/2)$, the
adjusted response can be either negative or greater than the binomial
total $m$. In such cases, standard implementations of maximum
likelihood are either unstable or report an error. This is because the
binomial log-likelihood is not necessarily concave when $y < 0$ or
$y > m$ for at least one observation, when a link function with
concave $\log\{G(\eta)\}$ and $\log\{1 - G(\eta)\}$ is used. 
Logit, probit, log-log and complementary
log-log are link functions of this kind.  See,
for example, \citet[][\S~5]{pratt:1981} 
for results and discussion on concavity of the
log-likelihood. 

Such issues with the use of repeated maximum-likelihood fits can be
avoided by noting that expression~(\ref{scores_jeffreys}) results if,
in the derivatives of the log-likelihood, $y$ and $m$ are replaced, respectively, by
their adjusted versions
\begin{equation}
  \label{pseudodatay}
  \tilde{y} = y + 2 a h (q - 1/2 + \pi c ) \quad \text{and} \quad  \tilde{m} = m + 2 a h c \, .
\end{equation}
Here $c$ is some arbitrarily chosen function of $\beta$.  The
following theorem identifies one function $c$ for which
$0 \le \tilde{y} \le \tilde{m}$.

\begin{theorem}
  \label{th:pseudodata}
  Let $I(A)$ be 1 if $A$ holds and 0 otherwise. If $c = 1 + (q - 1/2)\left\{ \pi - I(q \le 1/2) \right\}/\{\pi(1 - \pi)\}$, then $0 \le \tilde{y} \le \tilde{m}$.
\end{theorem}
The proof of Theorem~\ref{th:pseudodata} is given in the Supplementary
Material, which also provides pseudo-code (see Algorithm~S1)
and R code for Algorithm~\texttt{JeffreysMPL}, which implements
repeated maximum-likelihood fits to maximize the $l^\dagger(\beta; a)$
for any supplied $a$ and link function $G(\eta)$.

The variance-covariance matrix of the penalized likelihood estimator
can be obtained as $(R^\top R)^{-1}$, where $R$ is the upper
triangular matrix from the QR decomposition of $W(\beta)^{1/2} X$ at
the final iteration of the procedure. That decomposition is a
by-product of \texttt{JeffreysMPL}.

If, in addition to full rank $X$, we require that $X$ has a column of
ones and $g(\eta)$ is a unimodal density function, then it can be
shown that if the starting value of the parameter vector $\beta$ in
the repeated maximum-likelihood fits procedure has finite components,
then the values of $\beta$ computed in step P2 will also have finite
components at all repetitions. This is because, with a column of ones
in the full rank $X$, the adjusted responses and totals
in~(\ref{pseudodatay}) satisfy $0 < \tilde{y} < \tilde{m}$, and
hence maximum likelihood estimates with infinite components are not
possible. The strict inequalities $0 < \tilde{y} < \tilde{m}$ hold
because, under the aforementioned conditions, $w_i(\beta) > 0$ and
$X^\top W(\beta) X$ is positive definite for $\beta$ with finite
components. Then, \citet[Chapter~11, Theorem~4]{magnus+neudecker:1999}
on bounds for the Rayleigh quotient gives the inequality 
$h_i(\beta) \ge w_i(\beta) x_i^\top x_i \lambda(\beta) > 0$
$(1, \ldots, n)$, where $\lambda(\beta) > 0$ is the minimum eigenvalue
of $(X^\top W(\beta) X)^{-1}$.

The repeated maximum-likelihood fits procedure has the correct fixed
point even if, at step P2, full maximum-likelihood estimation is
replaced by a procedure that merely increases the log-likelihood, such
as a single step of iteratively reweighted least squares for the
adjusted responses and totals. \citet{firth:1992} suggested such a
scheme for logistic regressions with $a = 1/2$.  There is currently no
conclusive result on whether full maximum-likelihood iteration with a
reasonable stopping criterion is better or worse than, for example,
one step of iteratively reweighted least squares, in terms of
computational efficiency. A satisfactory starting value for the above
procedure is the maximum likelihood estimate of $\beta$, after adding
a small positive constant and twice that constant to the actual
binomial responses and totals, respectively.

Finally, for $a = 1/2$, repeated maximum-likelihood fits can be used
to compute the posterior normalizing constant when implementing the
importance sampling algorithm in \citet[\S~5]{chen+ibrahim+kim:2008}
for posterior sampling of the parameters of Bayesian binomial-response
generalized linear models with the Jeffreys prior.

Section~S3 of the Supplementary Material illustrates the
evolution of adjusted responses and totals through the iterations of
\texttt{JeffreysMPL}, for the first 6 games of Philadelphia 76ers in
Example~\ref{bt}. Section~S3 also computes the reduced-bias estimates
for a logistic regression model with $n = 1000$ binary responses and
$p = 200$ covariates, as considered in Figure~2(b) of the
supplementary information appendix of \citet{sur+candes:2019}, and
illustrates that such computation takes only a couple of seconds on a
standard laptop computer.

\section{Illustrations}

\begin{figure}[t!]
  \begin{center}
    \includegraphics[width = \textwidth]{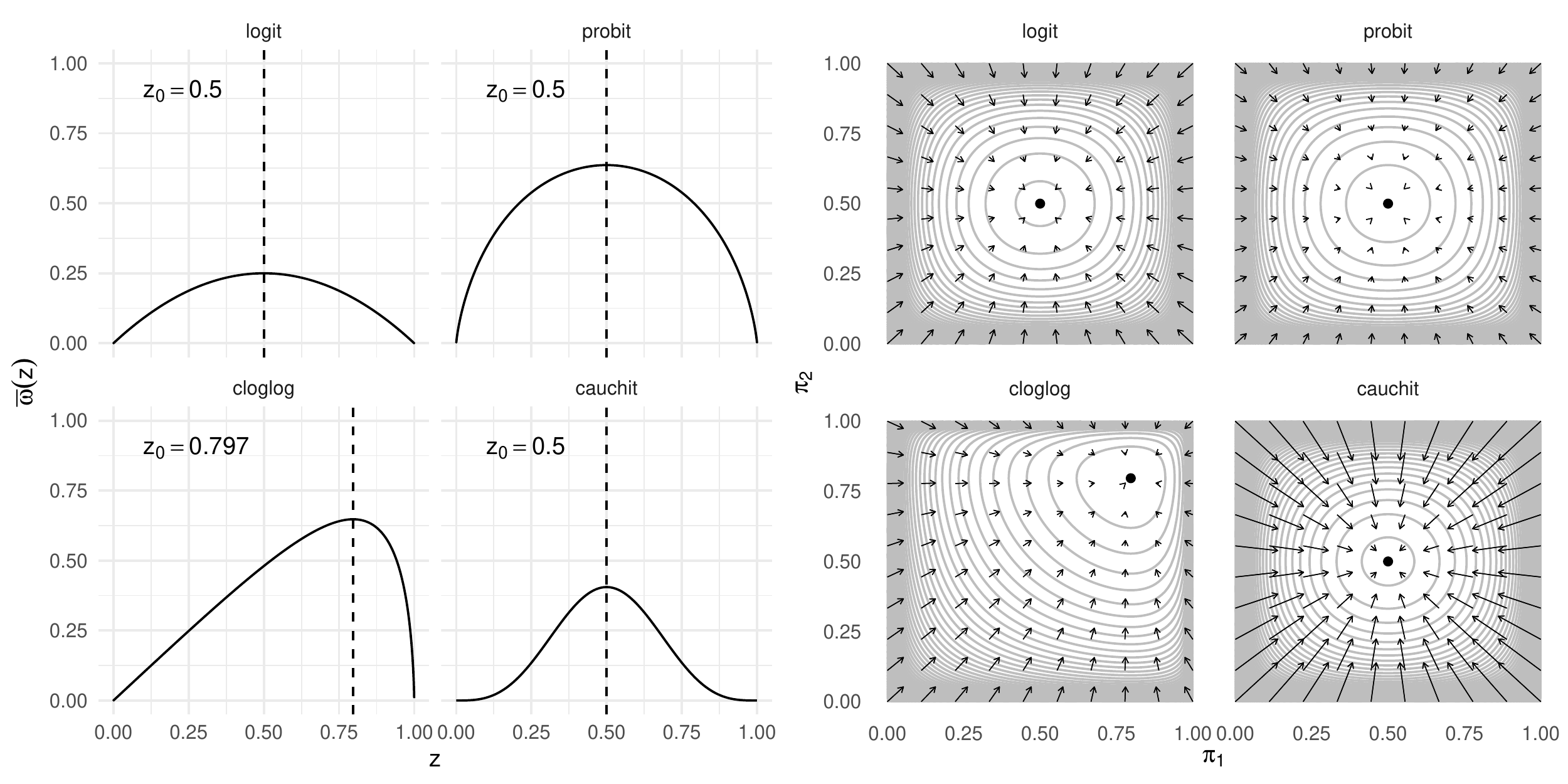}
  \end{center}
  \caption{Left: $\bar{\omega}(z)$ for various link functions. The
    dashed vertical line is at $z_0$. Right: Demonstration of how
    fitted probabilities from the penalized likelihood fit shrink
    relative to those of the maximum likelihood fit, from a complete
    enumeration of a saturated model where
    $\pi_i = G(\beta_1 + \beta_2 x_i)$ $(i = 1, 2)$, $x_1 = -1$,
    $x_2 = 1$ and $m_1 = 9$ and $m_2 = 9$. The arrows point from the
    estimated probabilities based on the maximum likelihood estimates
    to those based on the penalized likelihood estimates. The grey
    curves are the contours of $\log|X^\top \bar{W}(\pi)X|$. }
\label{shrinkage.plot}
\end{figure}

The left plot of Figure~\ref{shrinkage.plot} shows $\bar{\omega}(z)$
and $z_0$ for the various links. The plot for the log-log link is the
reflection of the one for the complementary log-log through $z =
0.5$. As is apparent, $\bar{\omega}(z)$ is concave for the logit,
probit and complementary-log-log links but not for the cauchit
link. The right plot of Figure~\ref{shrinkage.plot} visualizes the
shrinkage induced by the penalization by Jeffreys' invariant prior for
the logit, probit, complementary log-log and cauchit links. For each
link function, we obtain all possible fitted probabilities from a
complete enumeration of a saturated model with
$\pi_i = G(\beta_1 + \beta_2 x_i)$ $(i = 1, 2)$, where $x_1 = -1$,
$x_2 = 1$, $m_1 = 9$ and $m_2 = 9$.  The grey curves are the contours
of $\log|X^\top\bar{W}(\pi)X|$. An arrow is drawn from each pair of
estimated probabilities based on the maximum likelihood estimates to
the corresponding pair of estimated probabilities based on penalized
likelihood estimates, to demonstrate the induced shrinkage towards
$(z_0, z_0)^\top$ in accord to the results in Section~\ref{extensions}. Despite
the fact that $\bar\omega(z)$ is not concave for the cauchit link, the
fitted probabilities still shrink towards
$(z_0, z_0)^\top = (1/2, 1/2)^\top$. The plots in
Figure~\ref{shrinkage.plot} are invariant to the particular choice of
$x_1$ and $x_2$, as long as $x_1 \ne x_2$. For either maximum
likelihood or maximum penalized likelihood, if the estimates of
$\beta_1$ and $\beta_2$ are $b_1$ and $b_2$ for $x_1 = -1$ and
$x_2 = 1$, then the new estimates for any $x_1, x_2 \in \Re$ with
$x_1 \ne x_2$ are $b_1 - b_2 (x_1 + x_2) / (x_2 - x_1)$ and
$2 b_2 / (x_2 - x_1)$, respectively. Hence, the fitted probabilities
will be identical.

\begin{figure}[t!]
  \begin{center}
    \includegraphics[width = \textwidth]{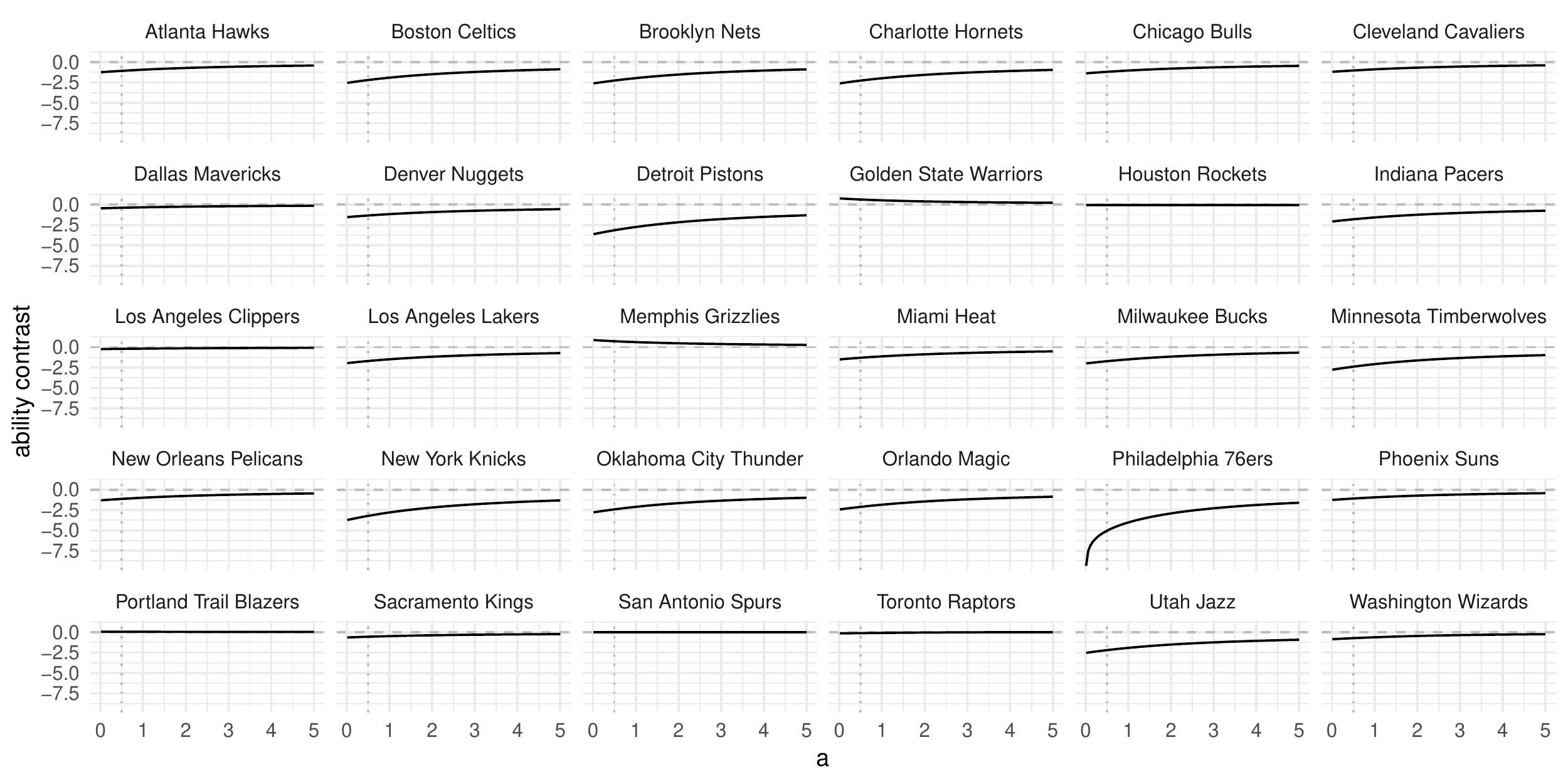}
  \end{center}
  \caption{Paths for the estimated ability contrasts from the
    maximization of~(\ref{penloglika}) for $a \in (0, 5]$. The dashed
    line is at zero and the dotted line is at $a = 0.5$, identifying
    the reduced-bias estimates on the paths.  }
    \label{btpaths}
  \end{figure}

Another illustration of finiteness and shrinkage follows from
Example~\ref{bt}. Figure~\ref{btpaths} shows the paths of the team
ability contrasts as $a$ varies from $0$ to $5$. The estimates are
obtained using \texttt{JeffreysMPL}, starting at the maximum
likelihood estimates of the ability contrasts after adding $0.01$ and
$0.02$ to the actual responses and totals, respectively. In accord
with the theoretical results in Section~\ref{finite}, the estimated
ability contrasts are finite for every $a > 0$; and, as expected from
the results in Section~\ref{shrinkage}, shrinkage towards
equiprobability becomes stronger as $a$ increases.

\section{Concluding remarks}

A recent stream of literature investigates the use of the coefficient
path defined by maximization of the penalized
log-likelihood~(\ref{penloglika}) for the prediction of rare events
through logistic
regression. \citet{elgmati+fiaccone+henderson+matthews:2015} study
that path for $a \in (0, 1/2]$, and propose to take $a$ to be around
$0.1$, in order to handle issues related to infinite estimates, and
they obtain predicted probabilities that are less biased than those
based on the reduced-bias estimates ($a = 0.5$). More recently,
\citet{puhr+heinze+nold+lusa+geroldinger:2017} proposed two new
methods for the prediction of rare events, and performed extensive
simulation studies to compare performance with various methods,
including maximum penalized likelihood with $a = 0.1$ and $a = 0.5$.

The coefficient path can be computed efficiently by using repeated
maximum-likelihood fits with ``warm'' starts. For a grid of values
$a_1 < \ldots < a_k$ with $a_j > 0$ $(j = 1, \ldots, k)$,
\texttt{JeffreysMPL} (Algorithm~S1 in the Supplementary Material) is
first applied with $a = a_1$ to get the maximum penalized likelihood
estimates $\beta^{(a_1)}$; then, \texttt{JeffreysMPL} is applied again
with $a = a_2$ with starting values $b = \beta^{(a_1)}$, and so on,
until $\beta^{(a_k)}$ has been computed. This process supplies
\texttt{JeffreysMPL} with the best available starting values, as the
algorithm walks through the grid. The finiteness of the components of
$\beta^{(a_1)}, \ldots, \beta^{(a_k)}$ and the shrinkage properties
described in Sections~\ref{shrinkage} and \ref{extensions} contribute
to the stability of the overall process. The properties of the
coefficient path for inference and prediction from binomial regression
models, and the development of general procedures for selecting $a$,
are interesting, open research topics.

\citet{kennepagui+salvan+sartori:2017} develop a method that can
reduce the median bias of the components of the maximum likelihood
estimator. According to the results therein, median bias reduction for
one-parameter logistic regression models is equivalent to
maximizing~(\ref{penloglika}) with $a = 1/6$. Hence, the results in
Section~\ref{logistic} also establish the finiteness of the estimate
from median bias reduction in one-parameter logistic regression, and
that the induced shrinkage to equiprobability will be less strong than
penalization by the Jeffreys prior.
\citet{kennepagui+salvan+sartori:2017} observed such properties
in numerical studies for $p > 1$.  When $p > 1$, though, 
median bias reduction is no longer equivalent to
maximizing~(\ref{penloglika}) with $a = 1/6$.

\section{Acknowledgments}
Ioannis Kosmidis and David Firth are supported by The Alan Turing
Institute under the EPSRC grant EP/N510129/1. David Firth was partly
supported also by EPSRC programme EP/K014463/1, \emph{Intractable
  Likelihood: New Challenges from Modern Applications}.

\section{Supplementary material}
The Supplementary Material is available for download at
\url{http://www.ikosmidis.com/files/finiteness-jeffreys-supplementary-v1.4.zip}
and includes: a document with proofs for Theorems 1, 2, 3 and
Corollary 1; Algorithm S1, and some additional numerical results; and
R code and data to reproduce all of the numerical work and graphs.

\bibliographystyle{biometrika}
\bibliography{finiteness_jeffreys}



\section*{Supplementary material}

# Supplementary material document for Jeffreys-prior penalty, finiteness and shrinkage in binomial-response generalized linear models

Ioannis Kosmidis
ioannis.kosmidis@warwick.ac.uk

and

David Firth
d.firth@warwick.ac.uk

Department of Statistics, University of Warwick
Coventry CV4 7AL, UK

and

The Alan Turing Institute
British Library, London NW1 2DB, UK

March 23, 2020

## S1  Supplementary material

All labels for the sections, equations, the table, the algorithm and the figure in the current document have been prefixed by "S" (e.g. Section S2, Table S1, Algorithm S1, etc).

The supplementary material for "Jeffreys-prior penalty, finiteness and shrinkage in binomial-response generalized linear models" provides i) the proofs of Theorem 1, Corollary 1, Theorem 2, Theorem 3 in the main text (see Section S2); ii) the pseudo-code for Algorithm JeffreysMPL (see Algorithm S1), which implements the repeated maximum-likelihood fits procedure of Section 4 in the main text to maximize the penalized log-likelihood $l^\dagger(\beta; a)$ in (4) for any supplied $a$ and link function $G(\eta)$; iii) illustrations using an R implementation of Algorithm S1; and iv) R scripts to reproduce all numerical results and figures in the main text and the current document. The current supplementary material document and the R scripts are available for download at http://www.ikosmidis.com/files/finiteness-jeffreys-supplementary-v1.3.zip.

In particular, the script nba-1415-case-study.R reproduces the numerical results in Example 1, Figure 1, and Figure 3 in the manuscript, and in Table S1 in the current document; the files nba-1415-functions.R and nba-1415-regular-season.csv provide the R functions and the data, respectively, for the case study in the main text; the script jeffreys-shrinkage.R reproduces Figure 2 in the main text; the file jeffreys-MPL.R provides an R implementation of Algorithm JeffreysMPL, and the file sur-candes-2019.R computes the timings for Algorithm JeffreysMPL that are reported in Section S3.3, and reproduces Figure 2b on page 11 of the supplementary information appendix of Sur and Candès (2019) (Figure S1 here).

All code has been tested and run in R version 3.6.3 (R Core Team, 2020) using the R packages brglm2 version 0.6.2 (Kosmidis, 2020a), enrichwith version 0.3.1 (Kosmidis, 2020b), ggplot2 version 3.2.1 (Wickham, 2016), qvcalc version 1.0.2 (Firth, 2020), rbenchmark version 1.0.0 (Kusnierczyk, 2012).

# S2 Proofs of Theorem 1, Corollary 1, Theorem 2, Theorem 3

## S2.1 Proof of Theorem 1

*Since $X$ has full rank, $|X^\top W(\beta)X|$ is not trivially zero for all $\beta \in \Re^p$. Let $R = \{1, \ldots, n\}$ and*

$$R_{(1)} = \{i : |\eta_i(\beta(r))| \to \infty \text{ as } r \to \infty; i \in R\}$$
$$R_{(2)} = \{i : \eta_i(\beta(r)) \to c_i, |c_i| < \infty \text{ as } r \to \infty; i \in R\}.$$

*Then $R = R_{(1)} \cup R_{(2)}$ and $R_{(1)} \cap R_{(2)} = \emptyset$, where $\emptyset$ is the empty set.*

*We first consider the case where $R_{(1)}$ and $R_{(2)}$ are non-empty. Then, $X$ and $W(\beta)$ can be partitioned as*

$$X = \begin{bmatrix} X_{(1)} \\ X_{(2)} \end{bmatrix} \quad \text{and} \quad W(\beta) = \begin{bmatrix} W_{(1)}(\beta) & 0 \\ 0 & W_{(2)}(\beta) \end{bmatrix},$$

*where $X_{(l)}$ has rows $x_i$ with $i \in R_{(l)}$ and $W_{(l)}(\beta) = \mathrm{diag}\{m_i\omega(\eta_i(\beta)), i \in R_{(l)}\}$ ($l = 1, 2$). We can then write*

$$X^\top W(\beta(r))X = X_{(1)}^\top W_{(1)}(\beta(r))X_{(1)} + X_{(2)}^\top W_{(2)}(\beta(r))X_{(2)}. \tag{S1}$$

*The limit of the first term in the right hand side of (S1), as $r \to \infty$, is zero because the limit of $\omega(\eta) = e^\eta/(1+e^\eta)^2$, as $\eta$ grows to infinity in absolute value, is zero. For the second term, $X_{(2)}$ is such that $X_{(2)}\beta(r) \to c$ as $r \to \infty$ where $c$ has all of its components finite. There must exist vectors $a \in \Re^p$ and $b \in \Re^p$ with finite components such that $a + br \to \beta_0$. So, $X_{(2)}a + X_{(2)}br \to c$ as $r \to \infty$ which is possible if and only if $X_{(2)}b = 0$. Hence $X_{(2)}$ has rank smaller than $p$ and $|X_{(2)}^\top W_{(2)}(\beta(r))X_{(2)}| = 0$ for all $r$. The result follows because $|X_{(2)}^\top W_{(2)}(\beta(r))X_{(2)}| = 0$ for all $r$ and $|X_{(1)}^\top W_{(1)}(\beta(r))X_{(1)}| \to 0$ as $r \to \infty$.*

*Because $X$ is of full rank and $X_{(2)}$ has rank smaller than $p$, $R_{(1)}$ cannot be empty. Hence, we only need to also examine the case that $R_{(2)}$ is empty and $R_{(1)}$ is not. In this case the same arguments as above give $|X^\top W(\beta(r))X| = |X_{(1)}^\top W_{(1)}(\beta(r))X_{(1)}| \to 0$ as $r \to \infty$.*

## S2.2 Proof of Corollary 1

*The binomial log-likelihood $l(\beta)$ in (2) is bounded above by zero. Hence, according to Theorem 1 and expression (3), $\tilde{l}(\beta(r)) \to -\infty$ as $\beta(r) \to \beta_0$. Such a setting for $\beta$ is always dominated, by a choice $b$ with finite components for which $\tilde{l}(b)$ takes a finite value. Hence, the maximizer of $\tilde{l}(\beta)$ must have finite components.*

## S2.3 Proof of Theorem 2

*The sum of $m$ independent Bernoulli distributions with probability $\pi$ is binomial with index $m$ and probability $\pi$. For this reason and without any loss of generality, the proof proceeds with $m_i = 1$ so that $w_i(\beta) = \omega(\eta_i(\beta))$ ($i = 1, \ldots, n$).*

*For proving i) decompose $X$ as $X = QR$, where $Q$ is a $n \times p$ matrix with orthonormal columns ($Q^\top Q = I_p$ where $I_p$ is the $p \times p$ identity matrix) and $R$ is a $p \times p$ non-singular matrix. This decomposition is always possible because $X$ has full rank by assumption. Then $|X^\top W(\beta)X| = |Q^\top W(\beta)Q||R|^2$. The functions $|X^\top W(\beta)X|$ and $|Q^\top W(\beta)Q|$ will have stationary points of the same kind and at the same values of $\beta$, because $|R|^2 > 0$ does not depend on $\beta$.*

*Denote the ordered set of quadratic weights as $w_{(1)}(\beta), \ldots, w_{(n)}(\beta)$ with $w_{(1)}(\beta) \leq \ldots \leq w_{(n)}(\beta)$. The Poincaré separation theorem (see, for example Magnus and Neudecker, 1999, Chapter 11, Theorem 10 for statement) and the positive definiteness of $W(\beta)$ for $\beta \in \Re^p$ imply that*

$$\prod_{t=1}^p w_{(t)}(\beta) \leq |Q^\top W(\beta)Q| \leq \prod_{t=1}^p w_{(n-p+t)}(\beta). \tag{S2}$$

*Note that $0 \leq \omega(\eta) \leq 1/4$, with the upper bound achieved when $\eta = 0$. It follows that $\prod_{t=1}^{p} w_{(t)}(\beta) \leq 1/4^p$ and $\prod_{t=1}^{p} w_{(n-t+t)}(\beta) \leq 1/4^p$ and that, at $\beta = 0$ inequalities (S2) become $1/4^p \leq |Q^\top W(0) Q| \leq 1/4^p$. The proof of i) concludes by noting that $|Q^\top W(0) Q| = 1/4^p$ which is the maximum value that $|Q^\top W(\beta) Q|$ can take.*

*For proving ii), note that $\bar{\omega}(z) = z(1-z)$ is concave. Hence, for $\theta \in (0,1)$, $\theta^* = 1 - \theta$ and any pair of n-vectors of probabilities $\pi$ and $\rho$ $\bar{\omega}(\theta \pi_i + \theta^* \rho_i) \geq \bar{\omega}(\theta \pi_i) + \bar{\omega}(\theta^* \rho_i)$ $(i = 1, \ldots, n)$. Lemma 1 can then be used to show that*

$$\left| X^\top \bar{W}(\theta \pi + \theta^* \rho) X \right| \geq \left| \theta X^\top \bar{W}(\pi) X + \theta^* X^\top \bar{W}(\rho) X \right|,$$

*The result in ii) follows from Magnus and Neudecker (1999, Chapter 11, Theorem 25) which can be used to show that*

$$\left| \theta X^\top \bar{W}(\pi) X + \theta^* X^\top \bar{W}(\rho) X \right| \geq \left| X^\top \bar{W}(\pi) X \right|^\theta \left| X^\top \bar{W}(\rho) X \right|^{\theta^*}.$$

**Lemma 1.** *If $A$ and $B$ are both diagonal $n \times n$ matrices with non-negative diagonal elements $\{a_r\}$ and $\{b_r\}$, respectively, and $a_r \geq b_r$, for every $r \in \{1, \ldots, n\}$, then, if $X$ is a $n \times p$ matrix, $|X^\top A X| \geq |X^\top B X|$.*

*Proof.* Since $A \geq B$, elementwise, $A = B + C$, where $C$ is a diagonal matrix with non-negative entries. Furthermore, $X^\top A X$, $X^\top B X$ and $X^\top C X$ are positive semidefinite, by the non-negativity of the diagonal elements of $A$, $B$ and $C$, respectively. Hence, by Magnus and Neudecker (1999, Chapter 11, Theorem 9), $\lambda_t \left( X^\top A X \right) \geq \lambda_t \left( X^\top B X \right)$ $(t = 1, \ldots, p)$, where $\lambda_t(D)$ denotes the $t$th eigenvalue of the matrix $D$. Since the determinant of a matrix is the product of its eigenvalues the result follows. □

### S2.4 Proof of Theorem 3

*For $c$ as in Theorem 3, the adjusted responses and totals in (6) have the form*

$$\tilde{y} = y + 2ah\pi \left\{ 1 + (q - 1/2) \frac{1 - I(q \leq 1/2)}{\pi(1 - \pi)} \right\}, \tag{S3}$$

$$\tilde{m} = m + 2ah \left\{ 1 + (q - 1/2) \frac{\pi - I(q \leq 1/2)}{\pi(1 - \pi)} \right\}.$$

*The result follows for any value of $q$, because $0 \leq h \leq 1$ and $0 \leq \pi \leq 1$.*

## S3 Algorithm `JeffreysMPL`

### S3.1 Details

Algorithm `JeffreysMPL` (see Algorithm S1) implements the repeated maximum-likelihood fits procedure of Section 4 in the main text, to maximize the penalized log-likelihood $l^\dagger(\beta; a)$ in (4) for any supplied $a$ and link function $G(\eta)$.

A satisfactory starting value $b$ for `JeffreysMPL` is the maximum likelihood estimate of $\beta$, after adding a small positive constant and twice that constant to the actual binomial responses and totals, respectively.

Step 22 of `JeffreysMPL` can be carried out using readily available maximum-likelihood implementations for binomial-response generalized linear models, such as the `glm` function in R (R Core Team, 2020) and the various implementations in the Python modules `statsmodels` (Seabold and Perktold, 2010) and `scikit-learn` (Pedregosa et al., 2011).

**Algorithm S1** Repeated maximum-likelihood fits for the maximization of $l^\dagger(\beta; a)$ in expression (4). The inputs y, m, X, a, G, g, gdash are $y = (y_1, \ldots, y_n)^\top$, $m = (m_1, \ldots, m_n)^\top$, $X$, $a$ in expression (4) of the main text, $G(\eta)$, $g(\eta)$ and $g'(\eta)$, respectively, and $||.||$ is the $L_2$ norm. The starting vector for $\beta$ is $b$ and $\epsilon$ is a small positive constant. ML is any maximum likelihood procedure.

1: **procedure** JEFFREYSMPL(y, m, X, a, G, g, gdash, ML, b, $\epsilon$)
2:     k ← 0
3:     n ← numberofrows(X)    ▷ Number of observations; must equal length(y) and length(m)
4:     p ← numberofcolumns(X)    ▷ Number of parameters; must equal length(b)
5:     eta, pi, d, dd ← vector(n)
6:     **for** i ∈ {1, 2, ..., n} **do**
7:       xi ← (X[i, 1], ..., X[i, p])
8:       eta[i] ← dotproduct(xi, b)    ▷ dot product of xi and b
9:       pi[i] ← G(eta[i])
10:       d[i] ← g(eta[i])
11:       dd[i] ← g'(eta[i])
12:     **end for**
13:     w ← d ∗ d/pi/(1 − pi)    ▷ elementwise operations
14:     q ← dd/w + pi    ▷ elementwise operations
15:     j ← I(q ≤ 1/2)    ▷ elementwise inequality and indicator function
16:     V ← diag($\sqrt{m[1] * w[1]}, \ldots, \sqrt{m[n] * w[n]}$) · X    ▷ · stands for matrix product
17:     Q ← orthogonal matrix from the QR decomposition of V
18:     h ← rowsums(Q ∗ Q)    ▷ sum the elements in each row of the elementwise product Q ∗ Q
19:     y_adj ← y + 2 ∗ a ∗ h ∗ pi ∗ (1 + (q − 1/2) ∗ (1 − j)/pi/(1 − pi)) ▷ elementwise operations
20:     m_adj ← m + 2 ∗ a ∗ h ∗ (1 + (q − 1/2) ∗ (pi − j)/pi/(1 − pi)) ▷ elementwise operations
21:     bp ← b
22:     b ← ML(y_adj, m_adj, X, G)    ▷ Maximum likelihood fit on adjusted data with link G
23:     **if** $||bp − b|| < \epsilon$ **then**
24:       **return** $b$
25:     **else**
26:       k ← k + 1
27:       Go to 6
28:     **end if**
29: **end procedure**

The estimated variance-covariance matrix of the penalized likelihood estimator can be obtained as $(R^\top R)^{-1}$, where $R$ is the upper triangular matrix from the QR decomposition of $W(\beta)^{1/2} X$ at the final iteration of the procedure. That decomposition is a by-product of step 17 of JeffreysMPL.

### S3.2 Illustration: evolution of adjusted responses and totals

Table S1 shows the values of the adjusted responses and totals for the first 6 games of Philadelphia 76ers in Example 1 of the main text, at the first 6 iterations of Algorithm JeffreysMPL, when computing the reduced-bias fit shown in Figure 1 (see script `nba-1415-case-study.R` for code to reproduce Table S1). The starting values (iteration 0) are the maximum likelihood estimates of the ability contrasts after adding 0.01 and 0.02 to the actual responses and totals, respectively.

Table S1: The adjusted responses (top) and totals (bottom) for the first 6 games of Philadelphia 76ers (P76) at the first 6 iterations of Algorithm S1, when computing the reduced-bias fit in Figure 1. Figures are shown in 3 decimal digits. The home team is mentioned first in column names. The actual response is 1 if the home team wins and 0 otherwise. The acronyms for the opponents are IP (Indiana Pacers), MB (Milwaukee Bucks), MH (Miami Heat), HR (Houston Rockets), OR (Orlando Magic) and CB (Chicago Bulls). The starting values are the maximum likelihood estimates of the ability contrasts after adding 0.01 and 0.02 to the actual responses and totals (iteration 0).

| Iteration | P76 vs IP | P76 vs MB | MH vs P76 | HR vs P76 | OM vs P76 | CB cs P76 |
|---|---|---|---|---|---|---|
| | | | Adjusted responses | | | |
| 0 | 0.010 | 0.010 | 1.010 | 1.010 | 1.010 | 1.010 |
| 1 | 0.039 | 0.036 | 1.045 | 1.012 | 1.104 | 1.039 |
| 2 | 0.045 | 0.042 | 1.054 | 1.017 | 1.110 | 1.048 |
| 3 | 0.046 | 0.043 | 1.055 | 1.017 | 1.111 | 1.049 |
| 4 | 0.046 | 0.043 | 1.055 | 1.018 | 1.111 | 1.049 |
| 5 | 0.046 | 0.043 | 1.055 | 1.018 | 1.111 | 1.049 |
| 6 | 0.046 | 0.043 | 1.055 | 1.018 | 1.111 | 1.049 |
| | | | Adjusted totals | | | |
| 0 | 1.020 | 1.020 | 1.020 | 1.020 | 1.020 | 1.020 |
| 1 | 1.114 | 1.105 | 1.067 | 1.018 | 1.158 | 1.059 |
| 2 | 1.128 | 1.120 | 1.081 | 1.025 | 1.170 | 1.073 |
| 3 | 1.130 | 1.122 | 1.083 | 1.026 | 1.171 | 1.075 |
| 4 | 1.131 | 1.122 | 1.084 | 1.026 | 1.172 | 1.075 |
| 5 | 1.131 | 1.122 | 1.084 | 1.026 | 1.172 | 1.075 |
| 6 | 1.131 | 1.122 | 1.084 | 1.026 | 1.172 | 1.075 |

### S3.3 Illustration: Algorithm `JeffreysMPL` for the high-dimensional logistic regression setting in Sur and Candès (2019)

The R implementation of `JeffreysMPL` in the supplementary material (see script `Jeffreys-MPL.R` for code) is used here to compute the reduced-bias estimates for a logistic regression model with $n = 1000$ binary responses and $p = 200$ covariates, as considered in Figure 2(b) of the supplementary information appendix of Sur and Candès (2019) (see script `sur-candes-2019.R` for reproducible code).

In particular, we construct a $1000 \times 200$ model matrix $X$ by simulating $200\,000$ independent random variables from a normal distribution with mean 0 and variance $10^{-3}$. Then, we simulate 1000 Bernoulli random variables from a logistic regression model with linear predictors $X\beta$, where $\beta_1 = \ldots = \beta_{25} = 10$, $\beta_{26} = \ldots = \beta_{50} = -10$, and $\beta_{50} = \ldots = \beta_{200} = 0$.

The R implementation of `JeffreysMPL` relies on a full maximum-likelihood iteration for step 22 using the R function `glm.fit`, and it takes approximately 2.73 seconds to converge to the reduced-bias estimates of the 200 parameters in 4 decimal places on a MacBook Pro laptop with 3.5GHz processor and 16GB of memory. The resulting maximum penalized likelihood estimates are shown in Figure S1 along with the corresponding maximum likelihood estimates.

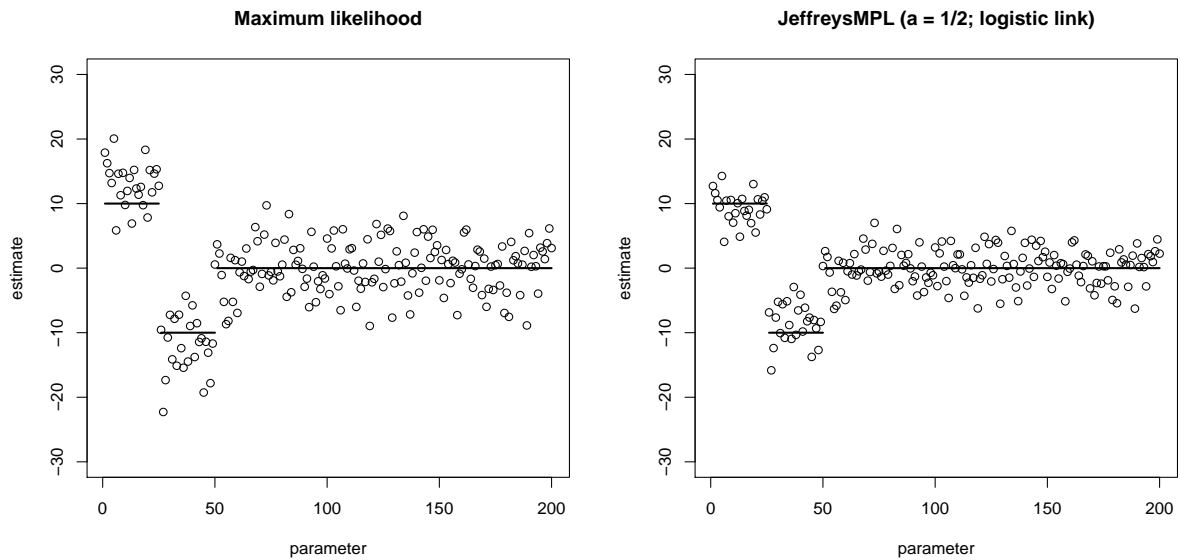

Figure S1: Maximum likelihood and reduced-bias estimates for the parameters of a logistic regression model on data simulated as in the supplementary information appendix of Sur and Candès (2019). In particular, we construct a $1000 \times 200$ model matrix $X$ by simulating $200\,000$ independent random variables from a normal distribution with mean 0 and variance $10^{-3}$. Then, we simulate 1000 Bernoulli random variables from a logistic regression model with linear predictors $X\beta$, where $\beta_1 = \ldots = \beta_{25} = 10$, $\beta_{26} = \ldots = \beta_{50} = -10$, and $\beta_{50} = \ldots = \beta_{200} = 0$. The maximum likelihood estimates are computed using the `glm` R function and the reduced-bias estimates are computed using Algorithm `JeffreysMPL`. The horizontal segments indicate the parameter values used for simulating the data.



\end{document}